\documentclass[12pt]{article}

\usepackage{amsmath,amssymb,amsbsy,amsfonts,amsthm,latexsym,
                        amsopn,amstext,amsxtra,euscript,amscd,color}
\usepackage{color,xcolor}
\usepackage{latexsym}
\usepackage{cite}
\usepackage{amsfonts}
\usepackage{amsfonts,mathrsfs}
\usepackage{graphicx}
\usepackage{psfrag}
\usepackage{subfigure}
\usepackage{url}
\usepackage{stfloats}
\usepackage{amsmath}

\newtheorem{theorem}{Theorem}
\newtheorem{lemma}{Lemma}

\newcommand{\quash}[1]{}

\begin{document}

\title{Linear complexity and trace representation of quaternary sequences over $\mathbb{Z}_4$ based on generalized cyclotomic classes modulo $pq$}

\author{Zhixiong Chen\\
 Key Laboratory of Applied Mathematics, Putian University, \\ Putian, Fujian
351100, P.R. China\\
ptczx@126.com
}

\maketitle

\begin{abstract}
We define a family of quaternary sequences over the residue class ring modulo $4$ of length $pq$, a product of two distinct odd primes, using the generalized cyclotomic classes modulo $pq$ and calculate the discrete Fourier transform (DFT) of the sequences.
The DFT helps us to determine the exact values of linear complexity and the trace representation of the sequences.

\end{abstract}

\section{Introduction}

Due to the applications of quaternary sequences in communication systems, radar and cryptography, see \cite{GG,HY}
it is of interest to design large families of quaternary sequences  over $\mathbb{Z}_4$.
There are many ways to define quaternary sequences.  A main method (for constructing quaternary sequences) is to use trace functions over Galois rings \cite{JHT,TU,US96,US98,US}. Second important way is to
use the inverse Gray mapping along two binary sequences \cite{EI15,KJKN,LNC,YK}. Using
cyclotomic and generalized cyclotomic classes is another important technique to define quaternary sequences \cite{TL,YT,YK}.
Most references concentrated on the correlation of the quaternary sequences. Linear complexity, as an important measure in cryptography, is also paid attention for certain quaternary sequences over $\mathbb{Z}_4$ \cite{EI13,EI14,US96}. However, it meets more difficulties due to the  phenomenon of zero divisors in $\mathbb{Z}_4$. In this manuscript, we will
define a family of quaternary sequences over $\mathbb{Z}_4$ by using the generalized cyclotomic classes modulo $pq$ and
investigate their linear complexities in terms of the discrete Fourier transform(DFT), from which we also derive the trace of the sequences.

 Let $m$ be a positive integer. We identify $\mathbb{Z}_m$,
the residue class ring modulo $m$, with the set $\{0,1,\cdots,m-1\}$ and
we denote by $\mathbb{Z}^*_m$  the unit group of $\mathbb{Z}_m$.

Let $p$ and $q$ be two distinct primes with $\gcd(p-1,q-1)=4$ and
$e=(p-1)(q-1)/4$. By the Chinese Remainder Theorem there exists a
common primitive root $g$ of both $p$ and $q$. There also exists an
integer $h$ satisfying
$$
h\equiv g\pmod{p},\,\  h\equiv 1\pmod{q}.
$$
Below we always fix the definitions of $g$ and $h$.  Since $g$ is a
primitive root of both $p$ and $q$, by the Chinese Remainder Theorem
again the multiplicative order of $g$
modulo $pq$ is $e$.

Define the generalized cyclotomic classes of order $4$ modulo $pq$ as
$$
D_i=\{g^sh^i \pmod {pq} : s=0,1,\ldots,e-1\},~~ 0\le i<4
$$
and we have
$$
\mathbb{Z}_{pq}^*=D_0\cup D_1 \cup D_2 \cup D_3.
$$
We note that $h^4\in D_0$, since otherwise, we write $h^4\equiv g^sh^i \bmod {pq}$ for some $0\le s<e$ and $1\le i<4$ and
get $g^{e-s}h^{4-i}=1 \in D_0$, a contradiction.

We also define
$$
P=\{p,2p,\ldots,(q-1)p\}, ~ Q=\{q,2q,\ldots,(p-1)q\}, ~R=\{0\}.
$$
Then we define the sequence $(e_u)$ over $\mathbb{Z}_4$ of length $pq$ by
\begin{equation}\label{quaternary}
e_u=\left\{
\begin{array}{ll}
2, & \mathrm{if}\,\ u\bmod {pq}\in Q\cup R,\\
0, & \mathrm{if}\,\ u\bmod {pq}\in P,\\
i, & \mathrm{if}\,\ u\bmod {pq}\in D_i, i=0,1,2,3.
\end{array}
\right.
\end{equation}
We remark that, except \cite{EI14}, most references mainly focused on binary sequences defined using   the generalized cyclotomic classes modulo $pq$, see e.g. \cite{CDR,DGS09,D97,D98,D98-3,DYX,YHX,YHX2}.

We organize this correspondence as follows. In Section 2, we calculate the Mattson-Solomon polynomial (see below for the definition) of $(e_u)$. We determine the linear complexity and present the trace representation of $(e_u)$ in terms of its   Mattson-Solomon polynomial in  Sections 3 and 4, respectively.

We conclude this section by introducing the notions of Galois rings of characteristic 4 and of the Mattson-Solomon polynomial of a quaternary sequence over $\mathbb{Z}_4$.

For a \emph{basic irreducible polynomial}  $f(X)\in \mathbb{Z}_4[X]$ of degree $r$, which means that $f(X)$ modulo 2 (i.e., the coefficients of $f(X)$ is reduced modulo 2) is irreducible over the finite field $\mathbb{F}_2$, the \emph{Galois ring} of characteristic 4 is defined as the residue class ring $\mathbb{Z}_4[X]/(f(X))$ of $4^r$ many elements and denoted by $GR(4,4^r)$. The group of units of $GR(4,4^r)$ satisfies
$$GR^*(4,4^r)=G_1\times G_2,$$
where $G_1$ is a cyclic group of order $2^r-1$ and $G_2$ is a group of order $2^r$. So we write $G_1=\langle \xi\rangle=\{\xi^i : 0\le i<2^r-1\}$ for some $\xi\in GR(4,4^r)$ of order $2^r-1$. Let
$$
\mathcal{T}=\{0\}\cup G_1=\left\{0,1,\xi,\xi^2,\ldots,\xi^{2^r-2}\right\},
$$
which is referred in the literature to as \emph{Teichmuller set}. Each element $\alpha\in GR(4,4^r)$  can be represented as
$$
\alpha := \alpha_1+2\alpha_2, ~~~ \alpha_1,\alpha_2\in\mathcal{T}.
$$
See \cite[Ch. 14]{W} for details on the theory of Galois rings.

For a quaternary sequence $(s_u)$ over $\mathbb{Z}_4$  of odd period $T$ (such that $T|(2^r-1)$ for some $r$), there exists a primitive $T$-th root
$\alpha\in GR(4,4^r)$ of unity such that
$$
s_u=\sum\limits_{0\le i<T}\rho_i\alpha^{iu},~~   0\le u<T,
$$
where $\rho_i$'s are given by
$$
\rho_i=\sum\limits_{0\le u<T}s_u\alpha^{-iu}, ~~   0\le i<T,
$$
see \cite{M,US}. In fact, this is an extension of the \emph{discrete Fourier transform} of binary sequences \cite[Ch. 6]{GG}.

If we write $G(X)= \sum\limits_{0\le i<T}\rho_iX^i \in  GR(4,4^r)[X]$, we have
$$
s_u=G(\alpha^u), ~~ u\ge 0.
$$
$G(x)$ is called the \emph{Mattson-Solomon polynomial} of $(s_u)$ in coding theory \cite{MS}.
Note that for a given $\alpha$, $G(X)$ is uniquely determined modulo $x^{T}-1$ since $T$ is co-prime to the characteristic of $GR(4,4^r)$.


\section{Mattson-Solomon polynomial of $(e_u)$}

It is easy to see that
$$
uD_i \triangleq \{uv \bmod{pq} : v\in D_i\}=D_{i+j}
$$
for $u\in D_j$. Here and hereafter the subscript of $D$ is performed modulo $4$, i.e., $D_{i+4}=D_i$ for all $0\le i<4$.

For a unit $\gamma\in GR^*(4,4^r)$, we denote by $\mathrm{ord}(\gamma)$ the order of $\gamma$, i.e., the least positive integer $n$ such that $\gamma^n=1$. From now on, we always suppose that the order of 2 modulo $pq$ is $\ell$, i.e., $\ell$ is the least number such that $2^\ell\equiv 1 \pmod {pq}$. So there exists a $\gamma \in GR(4,4^{\ell})$, the order of which is $\mathrm{ord}(\gamma)=pq$.

Define polynomials
$$
D_i(X)= \sum\limits_{u\in D_i}X^u \in
\mathbb{Z}_4[X]
$$
for $i=0,1,2,3$.

\begin{lemma}\label{lemma3}
Let $\gamma \in GR(4,4^{\ell})$ be a  primitive
$pq$-th root of unity, i.e., $\mathrm{ord}(\gamma)=pq$. We have

(1). $\gamma^{p}+\gamma^{2p}+\ldots+\gamma^{(q-1)p}=3$ or $1+\gamma^{p}+\gamma^{2p}+\ldots+\gamma^{(q-1)p}=0$.

(2). $\gamma^{q}+\gamma^{2q}+\ldots+\gamma^{(p-1)q}=3$ or $1+\gamma^{q}+\gamma^{2q}+\ldots+\gamma^{(p-1)q}=0$.

(3). $\sum\limits_{z\in \mathbb{Z}_{pq}^*}\gamma^z=1$ or $D_0(\gamma)+D_1(\gamma)+D_2(\gamma)+D_3(\gamma)=1$.

\end{lemma}
Proof. It is easy to check these results. We note that the calculations are performed in
the Galois ring $GR(4,4^{\ell})$ with characteristic four. ~\hfill $\square$

\begin{lemma}\label{lemma1}
Let $\gamma \in GR(4,4^{\ell})$ be a  primitive
$pq$-th root of unity. For $0\le i<4$, we have

(1). $D_i(1)=0$.

(2). $D_i(\gamma^{kq})=3(q-1)/4, ~~1\le k<p$.

(3). $D_i(\gamma^{kp})=3(p-1)/4, ~~1\le k<q$.
\end{lemma}
Proof. (1). Since each $D_i$ contains $(p-1)(q-1)/4$ many elements and $4|\frac{(p-1)(q-1)}{4}$, we have $D_i(1)=0$.

(2). We first compute
\begin{eqnarray*}
D_i \bmod p  & = &  \{(g^sh^i \bmod {pq}) \bmod p : s=0,1,\ldots,e-1\} \\
  & = & \{(g^{j+k(p-1)}h^i \bmod {pq}) \bmod p : 0\le j\le p-2, 0\le k<(q-1)/4\}\\
   & = & \{1,2,\ldots,p-1\},
\end{eqnarray*}
when $s$ ranges over $\{0,1,\ldots,e-1\}$, $g^sh^i \bmod p$ takes on each element of $\{1,2,\ldots,p-1\}$ $(q-1)/4$ times.
So for $1\le k<p$ we get by Lemma \ref{lemma3}(2)
$$D_i(\gamma^{kq})=\frac{q-1}{4}\sum\limits_{j\in \mathbb{Z}_{p}^*}\gamma^{jkq}=3(q-1)/4.$$

(3). The proof is similar to (2). ~ \hfill $\square$

\begin{lemma}\label{lemma2}
Let $0\le a<4$ be a fixed number.

(1). There are exactly $\frac{q-1}{4}$ many $w\in D_0$ such that $h^a+w \equiv 0 \pmod p$.

(2). There are exactly $\frac{p-1}{4}$ many $w\in D_0$ such that $h^a+w \equiv 0 \pmod q$.

(3). There is a unique $w\in D_0$ such that  $h^a+w \equiv 0 \pmod p$ and $h^a+w \equiv 0 \pmod q$ if and only if
$4|(\frac{p-1}{2}+a-\frac{q-1}{2})$.
\end{lemma}
Proof. (1). Let $w=g^x\in D_0$ for $0\le x<(p-1)(q-1)/4$. We have
$$
g^x\equiv -h^a\equiv g^{(p-1)/2+a} \pmod p,
$$
from which we derive $x\equiv (p-1)/2+a \pmod {p-1}$ and hence $x=k(p-1)+(p-1)/2+a$ for any $0\le k<(q-1)/4$.

One can prove (2) similarly.

For (3), we need to consider the equations
$$
\left\{
\begin{array}{ll}
x\equiv (p-1)/2+a & \pmod {p-1},\\
x\equiv (q-1)/2 & \pmod {q-1}.
\end{array}
\right.
$$
By \cite[Lemma 5]{D97}, $x$ exists iff
$4|(\frac{p-1}{2}+a-\frac{q-1}{2})$. If $x$ exists, $x$ is unique modulo $\frac{(p-1)(q-1)}{4}$.  ~ \hfill $\square$

We define $4$-tuples
$$
\mathcal{C}_i(X)=(D_{i}(X),D_{i+1}(X),D_{i+2}(X),D_{i+3}(X)), ~~ i=0,1,2,3.
$$
We  will calculate the inner product $\mathcal{C}_i(\beta)\cdot \mathcal{C}_j(\beta)$ for $0\le i,j<4$,
where $\beta \in GR(4,4^{\ell})$ is a  primitive
$pq$-th root of unity.

Since $\gcd(p-1,q-1)=4$, we see that $p$ and $q$ satisfy one of the following

-$p\equiv 1 \bmod 8$ and $q\equiv 5 \bmod 8$,

-$p\equiv 5 \bmod 8$ and $q\equiv 1 \bmod 8$,

-$p\equiv 5 \bmod 8$ and $q\equiv 5 \bmod 8$.

\begin{lemma}\label{inner-product}
Let $\beta \in GR(4,4^{\ell})$ be a  primitive
$pq$-th root of unity. For any fixed pair $0\le i,j<4$, we have
\begin{eqnarray*}
\mathcal{C}_i(\beta)\cdot \mathcal{C}_j(\beta) +\frac{q-1}{4}+\frac{p-1}{4}
=\left\{
\begin{array}{ll}
1, & \mathrm{if}\,\ i=j, \\
0, & \mathrm{otherwise},
\end{array}
\right.
\end{eqnarray*}
if $p\equiv 5 \bmod 8$ and $q\equiv 5 \bmod 8$, and
\begin{eqnarray*}
\mathcal{C}_i(\beta)\cdot \mathcal{C}_j(\beta) +\frac{q-1}{4}+\frac{p-1}{4}
=\left\{
\begin{array}{ll}
1, & \mathrm{if}\,\ (i,j)\in\{(2,0),(3,1),(0,2),(1,3)\}, \\
0, & \mathrm{otherwise},
\end{array}
\right.
\end{eqnarray*}
if $p\equiv 1 \bmod 8$ and $q\equiv 5 \bmod 8$ or $p\equiv 5 \bmod 8$ and $q\equiv 1 \bmod 8$.
\end{lemma}
Proof. Since $D_i=h^{i} D_0$ for all  $0\le i<4$, we calculate
\begin{eqnarray*}
\mathcal{C}_i(\beta)\cdot \mathcal{C}_j(\beta) & = & \sum\limits_{k=0}^{3}~\sum\limits_{u\in D_0}\beta^{uh^{i+k}}~\sum\limits_{v\in D_0}\beta^{vh^{j+k}}\\
          & = & \sum\limits_{k=0}^{3}~\sum\limits_{u\in D_0}\beta^{uh^{i+k}}~\sum\limits_{w\in D_0}\beta^{uwh^{j+k}}~~ ~~(\mathrm{we~ use~ } v=uw)\\
          & = & \sum\limits_{k=0}^{3}~\sum\limits_{u\in D_0}~\sum\limits_{w\in D_0}\beta^{uh^{j+k}(h^{i-j}+w)}\\
          & = & \sum\limits_{w\in D_0}~\sum\limits_{k=0}^{3}~~\sum\limits_{z\in D_{j+k}}\gamma_w^z ~~(\mathrm{we~ use~ } z=uh^{j+k}, \gamma_w=\beta^{h^{i-j}+w})\\
          & = & \sum\limits_{w\in D_0}~\sum\limits_{k=0}^{3}D_{k}(\gamma_w).
\end{eqnarray*}
Now we need to determine $\mathrm{ord}(\gamma_w)$, the order of $\gamma_w$ above for each $w\in D_0$. We note that $\mathrm{ord}(\gamma_w)|pq$ since $\beta$ is a primitive
$pq$-th root of unity, hence the possible values of $\mathrm{ord}(\gamma_w)$ are $1, p, q, pq$.

We first suppose that $p\equiv 5 \bmod 8$ and $q\equiv 5 \bmod 8$.

If $\mathrm{ord}(\gamma_w)= 1$, we find that $h^{i-j}+w \equiv 0 \pmod {pq}$. By Lemma \ref{lemma2}(3), there is a unique $w\in D_0$ satisfying this condition iff $4|(i-j)$ or $i=j$.

If $\mathrm{ord}(\gamma_w)= p$, we find that $h^{i-j}+w \equiv 0 \pmod {q}$ but $h^{i-j}+w \not\equiv 0\pmod {p}$.  By Lemma \ref{lemma2}(2), there are $\frac{p-1}{4}-1$ or $\frac{p-1}{4}$ many $w\in D_0$ satisfying this condition depending on whether $i=j$ or not.

Similarly, if $\mathrm{ord}(\gamma_w)= q$, we find that $h^{i-j}+w \equiv 0 \pmod {p}$ but $h^{i-j}+w \not\equiv 0\pmod {q}$.  By Lemma \ref{lemma2}(1), there are $\frac{q-1}{4}-1$ or $\frac{q-1}{4}$ many $w\in D_0$ satisfying this condition depending on whether $i=j$ or not.

So the number of $w\in D_0$  satisfying  $h^{i-j}+w \not\equiv 0\pmod {p}$ and $h^{i-j}+w \not\equiv 0\pmod {q}$, i.e., $\mathrm{ord}(\gamma_w)= pq$, is $\frac{(p-1)(q-1)}{4}-(\frac{q-1}{4}-1)-(\frac{p-1}{4}-1)-1$ or $\frac{(p-1)(q-1)}{4}-\frac{q-1}{4}-\frac{p-1}{4}$ depending on whether $i=j$ or not.

Putting everything together, we derive
\begin{eqnarray*}
\mathcal{C}_i(\beta)\cdot \mathcal{C}_j(\beta)  & = & \sum\limits_{\stackrel{w\in D_0}{\mathrm{ord}(\gamma_w)=1}}~\sum\limits_{k=0}^{3}D_{k}(\gamma_w)+\sum\limits_{\stackrel{w\in D_0}{\mathrm{ord}(\gamma_w)=p}}~\sum\limits_{k=0}^{3}D_{k}(\gamma_w)+\\
&&~~ \sum\limits_{\stackrel{w\in D_0}{\mathrm{ord}(\gamma_w)=q}}~\sum\limits_{k=0}^{3}D_{k}(\gamma_w)+\sum\limits_{\stackrel{w\in D_0}{\mathrm{ord}(\gamma_w)=pq}}~\sum\limits_{k=0}^{3}D_{k}(\gamma_w)
\end{eqnarray*}
and hence
\begin{eqnarray*}
\mathcal{C}_i(\beta)\cdot \mathcal{C}_j(\beta)  +\frac{q-1}{4}+\frac{p-1}{4}
=\left\{
\begin{array}{ll}
1, & \mathrm{if}\,\ i=j, \\
0, & \mathrm{otherwise}.
\end{array}
\right.
\end{eqnarray*}

For the case of $p\equiv 1 \bmod 8$ and $q\equiv 5 \bmod 8$ or $p\equiv 5 \bmod 8$ and $q\equiv 1 \bmod 8$,
one can derive the following result in a similar way
\begin{eqnarray*}
\mathcal{C}_i(\beta)\cdot \mathcal{C}_j(\beta)  +\frac{q-1}{4}+\frac{p-1}{4}
=\left\{
\begin{array}{ll}
1, & \mathrm{if}\,\ (i,j)\in\{(2,0),(3,1),(0,2),(1,3)\}, \\
0, & \mathrm{otherwise}.
\end{array}
\right.
\end{eqnarray*}
We note that in this case $h^{i-j}+w \equiv 0 \pmod {pq}$ has a (unique) solution $w\in D_0$ iff $4|(2+i-j)$ or $(i,j)\in\{(2,0),(3,1),(0,2),(1,3)\}$ by Lemma \ref{lemma2}(3).   ~\hfill $\square$

Now we present our main results.
\begin{theorem}\label{def-poly-caseone}
Let $\beta \in GR(4,4^{\ell})$ be a  primitive
$pq$-th root of unity. If $p\equiv 5 \bmod 8$ and $q\equiv 5 \bmod 8$, then the Mattson-Solomon polynomial $G(X)$ (corresponding to $\beta$) of the quaternary sequence $(e_u)$ over $\mathbb{Z}_4$ defined in (\ref{quaternary})
is
$$
G(X)=2\sum\limits_{j=0}^{p-1}X^{jq}+\sum\limits_{k=0}^{3}(\rho-k) D_k(X),
$$
where  $\rho=D_1(\beta)+ 2D_2(\beta)+ 3D_3(\beta)$.
\end{theorem}
Proof. By Lemma \ref{inner-product}, One can check that the defining polynomial $G(X)$
of $(e_u)$ is
\begin{eqnarray*}
G(X)& =&
2\sum\limits_{j=0}^{p-1}X^{jq}+\left( \mathcal{C}_1(\beta)\cdot \mathcal{C}_0(X) +\frac{q-1}{4}+\frac{p-1}{4}\right)\\
&   & \qquad + 2\left(\mathcal{C}_2(\beta)\cdot \mathcal{C}_0(X)+\frac{q-1}{4}+\frac{p-1}{4}\right)\\
&   & \qquad\quad +3\left(\mathcal{C}_3(\beta)\cdot \mathcal{C}_0(X)+\frac{q-1}{4}+\frac{p-1}{4}\right)\\
& =& 2\sum\limits_{j=0}^{p-1}X^{jq}+ \sum\limits_{i=1}^{3}
i\mathcal{C}_i(\beta)\cdot \mathcal{C}_0(X).
\end{eqnarray*}
In fact, for $u=0$, since $\mathcal{C}_0(1)=(0,0,0,0)$ by Lemma \ref{lemma1}(1), we have
$$
G(\beta^{0})=2\sum\limits_{j=0}^{p-1}1+\sum\limits_{i=1}^{3}
i\mathcal{C}_i(\beta)\cdot \mathcal{C}_0(\beta^{0})
=2p+0=2=e_{0}.
$$

For $u= kp$ with $1\le k<q$, we have
$$
\mathcal{C}_i(\beta)\cdot \mathcal{C}_0(\beta^{kp}) =\frac{3(p-1)}{4}(D_0(\beta)+D_1(\beta)+D_2(\beta)+D_3(\beta))=\frac{3(p-1)}{4}
$$
since $\mathcal{C}_0(\beta^{kp})=(\frac{3(p-1)}{4},\frac{3(p-1)}{4},\frac{3(p-1)}{4},\frac{3(p-1)}{4})$ by Lemma \ref{lemma1}(3) and hence
$$
G(\beta^{kp})=2\sum\limits_{j=0}^{p-1}1+\sum\limits_{i=1}^{3}
i\mathcal{C}_i(\beta)\cdot \mathcal{C}_0(\beta^{kp})
=2p+\frac{3(p-1)}{2}=2+2=0=e_{kp}.
$$

Similarly for $u= kq$ with $1\le k<p$, due to $\mathcal{C}_0(\beta^{kq})=(\frac{3(q-1)}{4},\frac{3(q-1)}{4},\frac{3(q-1)}{4},\frac{3(q-1)}{4})$ we have by Lemma \ref{lemma3}(2)
$$
G(\beta^{kq})=2\sum\limits_{j=0}^{p-1}\beta^{jkq^2}+\sum\limits_{i=1}^{3}
i\mathcal{C}_i(\beta)\cdot \mathcal{C}_0(\beta^{kq})
=0+\frac{3(q-1)}{2}=2=e_{kq}.
$$

For $u\in D_{k}$ with $0\le k<4$, we have by Lemmas \ref{lemma3}(2) and \ref{inner-product}
 \begin{eqnarray*}
G(\beta^u)& =&2\sum\limits_{j=0}^{p-1}\beta^{ujq}+\sum\limits_{i=1}^{3}
i\mathcal{C}_i(\beta)\cdot \mathcal{C}_0(\beta^{u})  \\
                  & = &0+\sum\limits_{i=1}^{3}
i\mathcal{C}_i(\beta)\cdot \mathcal{C}_k(\beta)  \\
                  & = &\sum\limits_{i=1}^{3}i\left( \mathcal{C}_i(\beta)\cdot \mathcal{C}_k(\beta) +\frac{q-1}{4}+\frac{p-1}{4}\right)\\
              & = & k=e_u.
\end{eqnarray*}
Hence we get $e_u=G(\beta^u)$ for all $u\ge 0$.

On the other hand, we re-write $G(X)$ as
\begin{eqnarray*}
G(X)& = & 2\sum\limits_{j=0}^{p-1}X^{jq}+ (D_1(\beta)+ 2D_2(\beta)+ 3D_3(\beta))D_0(X)\\
    &   & \qquad +(3D_0(\beta)+ D_2(\beta)+ 2D_3(\beta))D_1(X)\\
    &   & \qquad\quad +(2D_0(\beta)+ 3D_1(\beta)+ D_3(\beta))D_2(X)\\
    &   & \qquad\qquad +(D_0(\beta)+ 2D_1(\beta)+ 3D_2(\beta))D_3(X).
\end{eqnarray*}
Since $\rho=D_1(\beta)+ 2D_2(\beta)+ 3D_3(\beta)$, from Lemma \ref{lemma3}(3) we get
$$
\left\{
\begin{array}{l}
3D_0(\beta)+ D_2(\beta)+ 2D_3(\beta)=D_1(\beta)+ 2D_2(\beta)+ 3D_3(\beta)-\sum\limits_{i=0}^{3}D_i(\beta)=\rho-1, \\
2D_0(\beta)+ 3D_1(\beta)+ D_3(\beta)=3D_0(\beta)+ D_2(\beta)+ 2D_3(\beta)-\sum\limits_{i=0}^{3}D_i(\beta)=\rho-2,\\
D_0(\beta)+ 2D_1(\beta)+ 3D_2(\beta)=2D_0(\beta)+ 3D_1(\beta)+ D_3(\beta)-\sum\limits_{i=0}^{3}D_i(\beta)=\rho-3.
\end{array}
\right.
$$
This completes the proof. ~\hfill $\square$

Using a method similar to the one in proving Theorem \ref{def-poly-caseone}, one can obtain the Mattson-Solomon polynomial $G(X)$
of $(e_u)$ if $p\equiv 1 \bmod 8$ and $q\equiv 5 \bmod 8$ or $p\equiv 5 \bmod 8$ and $q\equiv 1 \bmod 8$.

\begin{theorem}\label{def-poly-casetwo}
Let $\beta \in GR(4,4^{\ell})$ be a  primitive
$pq$-th root of unity. If $p\equiv 1 \bmod 8$ and $q\equiv 5 \bmod 8$, then the Mattson-Solomon polynomial $G(X)$ (corresponding to $\beta$) of the quaternary sequence $(e_u)$ over $\mathbb{Z}_4$ defined in (\ref{quaternary})
is
$$
G(X)=2\sum\limits_{j=0}^{p-1}X^{jq}+2\sum\limits_{j=1}^{q-1}X^{jp}+\sum\limits_{k=0}^{3}(\rho+2-k) D_k(X),
$$
where  $\rho=D_1(\beta)+ 2D_2(\beta)+ 3D_3(\beta)$.
\end{theorem}

\begin{theorem}\label{def-poly-casethree}
Let $\beta \in GR(4,4^{\ell})$ be a  primitive
$pq$-th root of unity. If $p\equiv 5 \bmod 8$ and $q\equiv 1 \bmod 8$, then the Mattson-Solomon polynomial $G(X)$ (corresponding to $\beta$) of the quaternary sequence $(e_u)$ over $\mathbb{Z}_4$ defined in (\ref{quaternary})
is
$$
G(X)= 2+\sum\limits_{k=0}^{3}(\rho+2-k) D_k(X),
$$
where  $\rho=D_1(\beta)+ 2D_2(\beta)+ 3D_3(\beta)$.
\end{theorem}

\section{Linear complexity}

We recall that the
\emph{linear complexity} $LC((s_u))$  of a quaternary sequence $(s_u)$ over $\mathbb{Z}_4$ with period $T$ is the least order $L$ of a linear
recurrence relation over $\mathbb{Z}_4$
$$
s_{u+L} + c_{1}s_{u+L-1} +\ldots +c_{L-1}s_{u+1}+ c_Ls_u=0\quad
\mathrm{for}\,\ u \geq 0,
$$
which is satisfied by $(s_u)$ and where $c_1, c_2, \ldots,
c_{L}\in \mathbb{Z}_4$. The \emph{connection polynomial} is $C(X)$ given by $1+c_1X+\ldots+c_LX^L$.
Let
$S(X)=s_0+s_1X+\ldots+s_{T-1}X^{T-1}\in \mathbb{Z}_4[X]$ be the \emph{generating polynomial} of $(s_u)$. Then
an LFSR with a connection polynomial $C(X)$ generates  $(s_u)$, if and only if,
$$
S(X)C(X) \equiv 0 \pmod {X^T-1}.
$$
That is,
$$LC((s_u))=\min\{\deg(C(X)) : S(X)C(X) \equiv 0 \pmod {X^T-1}\}.$$
If $T$ is odd, Udaya and Siddiqi proved in \cite[Theorem 4]{US} that
 the linear complexity $LC((s_u))$ equals the number of nonzero coefficients of the Mattson-Solomon polynomial $G(X)$ of $(s_u)$.

\begin{theorem}\label{linearcomp}
The linear complexity $LC((e_u))$ over $\mathbb{Z}_4$  of the quaternary sequence $(e_u)$ defined in (\ref{quaternary})
is
$$
LC((e_u))= p+(p-1)(q-1) ~~~ \mathrm{or} ~~~  p+(p-1)(q-1)-\frac{(p-1)(q-1)}{4},
$$
if $p\equiv 5 \bmod 8$ and $q\equiv 5 \bmod 8$; and
$$
LC((e_u))= p+q-1+(p-1)(q-1) ~~~ \mathrm{or} ~~~  p+q-1+(p-1)(q-1)-\frac{(p-1)(q-1)}{4},
$$
If $p\equiv 1 \bmod 8$ and $q\equiv 5 \bmod 8$; and
$$
LC((e_u))= 1+(p-1)(q-1) ~~~ \mathrm{or} ~~~  1+(p-1)(q-1)-\frac{(p-1)(q-1)}{4},
$$
if $p\equiv 5 \bmod 8$ and $q\equiv 1 \bmod 8$.
\end{theorem}
Proof. We only need to consider the coefficients of $G(X)$ in Theorems \ref{def-poly-caseone}, \ref{def-poly-casetwo} and \ref{def-poly-casethree} are whether zero or not. We find that at most one element of $\rho,\rho-1,\rho-2,\rho-3$ is zero. So the desired result follows by \cite[Theorem 4]{US}. ~\hfill $\square$

\section{Trace representation}

Let $\phi: GR(4,4^{r}) \rightarrow GR(4,4^{r})$ be the \emph{Frobenius automorphism} defined by
$$
\phi : \alpha_1+2\alpha_2 \mapsto \alpha_1^2+2\alpha_2^2, ~~~ \alpha_1,\alpha_2\in\mathcal{T},
$$
where $\mathcal{T}$ is  the Teichmuller set of $GR(4,4^{r})$. Let $\phi^0=1$ be the identity map of $GR(4,4^{r})$ and $\phi^j=\phi^{j-1}\circ \phi$ for $j\ge 2$. Then the \emph{trace function} $\mathrm{TR}_{1}^{r}(-)$ from $ GR(4,4^{r})$ to $\mathbb{Z}_4$ is
defined by
$$
\mathrm{TR}_{1}^{r}(\alpha)=\phi^0(\alpha) +\phi(\alpha)+\ldots+\phi^{r-1}(\alpha), ~~~ \alpha\in GR(4,4^{r}).
$$
If $s|r$, the \emph{generalized Frobenius automorphism} of $GR(4,4^{r})$ over $GR(4,4^{s})$ is defined by
$$
\Phi_s  : \alpha_1+2\alpha_2 \mapsto \alpha_1^{2^s}+2\alpha_2^{2^s}, ~~~ \alpha_1,\alpha_2\in\mathcal{T},
$$
then the \emph{generalized trace function} $\mathrm{TR}_{s}^{r}(-)$ from $ GR(4,4^{r})$ to $ GR(4,4^{s})$ is
defined by
$$
\mathrm{TR}_{s}^{r}(\alpha)=\Phi_s^0(\alpha) +\Phi_s(\alpha)+\ldots+\Phi_s^{r/s-1}(\alpha), ~~~ \alpha\in GR(4,4^{r}).
$$
We note that the order of $\Phi_s$ is $r/s$, i.e., $\Phi_s^{r/s}=1$. In particular, for any $\alpha_1\in\mathcal{T}$ we have
$$
\mathrm{TR}_{s}^{r}(\alpha_1)=\alpha_1 +\alpha_1^{2^s}+\ldots+\alpha_1^{2^{s(r/s-1)}}.
$$
For more details on trace functions over Galois rings, we refer the reader to \cite{W}. The trace functions play an important role in sequences design \cite{GG,US}.

\begin{lemma}\label{2inDj}
(1). $2\in D_0\cup D_2$ if and only if $p\equiv 5 \bmod 8$ and $q\equiv 5 \bmod 8$.

(2).  $2\in D_1\cup D_3$ if and only if $p\equiv 1 \bmod 8$ and $q\equiv 5 \bmod 8$, or $p\equiv 5 \bmod 8$ and $q\equiv 1 \bmod 8$.
\end{lemma}
Proof.  Let $2\equiv g^sh^i\bmod {pq}$ for some $0\le s<\frac{(p-1)(q-1)}{4}$ and $0\le i<4$. If $i$ is even,  we see that
both $s$ and $s+i$ are odd or even. Hence $2$ is quadratic non-residue or quadratic residue modulo $p$ and $q$ respectively,
which means that $p\equiv 5 \bmod 8$ and $q\equiv 5 \bmod 8$ under the assumption of $\gcd(p-1,q-1)=4$. Conversely,
if  $p\equiv 5 \bmod 8$ and $q\equiv 5 \bmod 8$, $2$ is quadratic non-residue modulo $p$ and $q$, respectively.
From $2\equiv g^sh^i\equiv g^{s+i}\bmod {p}$  and $2\equiv g^sh^i\equiv g^s \bmod {q}$, we see that both $s$ and $s+i$ are odd, which indicates $i$ is even. We prove (1).

For odd $i$, we can prove (2) in a similar way.  ~ \hfill $\square$

\begin{lemma}\label{ell-2-4}
Let $\ell$ be the order of 2 modulo $pq$. We have

(1). $2|\ell$ if $2\in D_2$.

(2). $4|\ell$ if $2\in D_1\cup D_3$.
\end{lemma}
Proof. Let $2\equiv g^sh^i\bmod {pq}$ for some $0\le s<\frac{(p-1)(q-1)}{4}$ and $1\le i\le 3$.  Then from $1\equiv 2^{\ell}\equiv g^{s\ell}h^{i\ell}\bmod {pq}\in D_0$, we see that $h^{i\ell}\in D_0$ and hence $4|i\ell$, which implies the desired results. ~ \hfill $\square$

\begin{theorem}\label{trace-caseone}
With notations $g,h,e$ as in Section 1. Let $\ell$ be the order of 2 modulo $pq$ and $\ell_p$  the order of 2 modulo $p$. Let $\beta \in GR(4,4^{\ell})$ be a  primitive
$pq$-th root of unity. If $p\equiv 5 \bmod 8$ and $q\equiv 5 \bmod 8$, then the trace representation of the quaternary sequence $(e_u)$ over $\mathbb{Z}_4$ defined in (\ref{quaternary})
is
$$
e_u=2+2\sum\limits_{i=0}^{\frac{p-1}{\ell_p}-1}\mathrm{TR}_{1}^{\ell_p}(\beta^{ug^iq})+\sum\limits_{j=0}^{3}(\rho-j)~ \sum\limits_{i=0}^{\frac{e}{\ell}-1}\mathrm{TR}_{1}^{\ell}(\beta^{ug^ih^j})
$$
if $2\in D_0$, and
$$
e_u=2+2\sum\limits_{i=0}^{\frac{p-1}{\ell_p}-1}\mathrm{TR}_{1}^{\ell_p}(\beta^{ug^iq})+\sum\limits_{j=0}^{3}(\rho-j)~ \sum\limits_{i=0}^{\frac{2e}{\ell}-1}\mathrm{TR}_{2}^{\ell}(\beta^{ug^ih^j})
$$
if $2\in D_2$,
where  $\rho=D_1(\beta)+ 2D_2(\beta)+ 3D_3(\beta)$.
\end{theorem}
Proof. We only need to describe $D_i(X)$ and $\sum\limits_{j=1}^{p-1}X^{jq}$ using (generalized) trace functions over Galois rings.

First let
$$
U_p=\{2^{j} \pmod {p} : 0\le j< \ell_p\}\subseteq \mathbb{Z}_p^*,
$$
then $\mathbb{Z}_p^*$ is divided into $(p-1)/\ell_p$ many disjoint subsets
$$
U_p, ~ gU_p, \ldots, g^{(p-1)/\ell_p-1}U_p.
$$
So we have in $\mathbb{Z}_4[X]$
$$
U_p(X)=\sum\limits_{u\in U_p} X^u=X +X^2+\ldots+X^{2^{\ell_p-1}}\pmod {X^p-1}
$$
and
$$
\sum\limits_{j=1}^{p-1}X^{jq}=\sum\limits_{i=0}^{\frac{p-1}{\ell_p}-1}U_p\left(X^{g^{i}q}\right)\pmod {X^p-1}.
$$

Now if $2\in D_0$, let
$$
U=\{2^{j} \pmod {pq} : 0\le j< \ell\}\subseteq D_0.
$$
then $D_0$ is divided into $e/\ell$ many disjoint subsets
$$
U, ~ gU, \ldots, g^{e/\ell-1}U.
$$
Applying
$$
U(X)=\sum\limits_{u\in U} X^u=X +X^2+\ldots+X^{2^{\ell-1}} \pmod {X^{pq}-1}
$$
we derive
$$
D_0(X)=\sum\limits_{j=0}^{e/\ell-1}U\left(X^{g^{j}}\right) \pmod {X^{pq}-1}
$$
and
$$
D_i(X)=\sum\limits_{j=0}^{e/\ell-1}U\left(X^{g^{j}h^i}\right)\pmod {X^{pq}-1}
$$
for $1\le i<4$.

If $2\in D_2$, we have $4\in D_0$ and the order of 4 modulo $pq$ is $\ell/2$ by Lemma \ref{ell-2-4}(1). So let
$$
V=\{4^{j} \pmod {pq} : 0\le j< \ell/2\}\subseteq D_0.
$$
Then $D_0$ is divided into $2e/\ell$ many disjoint subsets
$$
V, ~ gV, \ldots, g^{2e/\ell-1}V.
$$
We can use
$$
V(X)=\sum\limits_{u\in V} X^u=X +X^4+\ldots+X^{4^{\ell/2-1}}\pmod {X^{pq}-1},
$$
to describe
$$
D_i(X)=\sum\limits_{j=0}^{2e/\ell-1}V\left(X^{g^{j}h^i}\right) \pmod {X^{pq}-1}
$$
for $0\le i<4$.

So for $\beta \in GR(4,4^{\ell})$ of order $pq$, we use the trace representations
$$
U_p(\beta^q)=\mathrm{TR}^{\ell_p}_{1}(\beta^q), ~~   U(\beta)=\mathrm{TR}_{1}^{\ell}(\beta),   ~~ V(\beta)=\mathrm{TR}_{2}^{\ell}(\beta)
$$
to complete the proof. ~\hfill $\square$

\begin{theorem}\label{trace-casetwo}
With notations $g,h,e$ as in Section 1. Let $\ell$ be the order of 2 modulo $pq$, $\ell_p$  the order of 2 modulo $p$ and $\ell_q$  the order of 2 modulo $q$. Let $\beta \in GR(4,4^{\ell})$ be a  primitive
$pq$-th root of unity. Then the trace representation of the quaternary sequence $(e_u)$ over $\mathbb{Z}_4$ defined in (\ref{quaternary})
is
$$
e_u=2+2\sum\limits_{i=0}^{\frac{p-1}{\ell_p}-1}\mathrm{TR}_{1}^{\ell_p}(\beta^{ug^iq})
+2\sum\limits_{i=0}^{\frac{q-1}{\ell_q}-1}\mathrm{TR}_{1}^{\ell_q}(\beta^{ug^ip})
+\sum\limits_{j=0}^{3}(\rho+2-j)~ \sum\limits_{i=0}^{\frac{4e}{\ell}-1}\mathrm{TR}_{4}^{\ell}(\beta^{ug^ih^j}),
$$
if $p\equiv 1 \bmod 8$ and $q\equiv 5 \bmod 8$; while
$$
e_u=2+\sum\limits_{j=0}^{3}(\rho+2-j)~\sum\limits_{i=0}^{\frac{4e}{\ell}-1}\mathrm{TR}_{4}^{\ell}(\beta^{ug^ih^j}),
$$
if $p\equiv 5 \bmod 8$ and $q\equiv 1 \bmod 8$,
where  $\rho=D_1(\beta)+ 2D_2(\beta)+ 3D_3(\beta)$.
\end{theorem}
Proof. The proof is similar to that of Theorem \ref{trace-caseone}, here we only present some sketch.
since $2\in D_1\cup D_3$ by  Lemma \ref{2inDj}(2), we have $16\in D_0$ and the order of 16 modulo $pq$ is $\ell/4$ by Lemma \ref{ell-2-4}(2). Let
$$
W=\{16^{j} \pmod {pq} : 0\le j< \ell/4\}\subseteq D_0.
$$
Then for $\beta \in GR(4,4^{\ell})$ of order $pq$, we use the trace representation
$$
W(\beta)=\mathrm{TR}_{4}^{\ell}(\beta)
$$
to complete the proof. ~\hfill $\square$

\section{Conclusions}

Determining linear complexity of  quaternary sequences over $\mathbb{Z}_4$ is a bottleneck problem due to the zero divisors of
$\mathbb{Z}_4$. It is interesting to develop a way to solve this problem.
In this work, we define a special family of quaternary sequences over $\mathbb{Z}_4$ using the generalized cyclotomic classes modulo $pq$ and determine the linear complexities by computing their discrete Fourier transform.
We also give the trace representation of the sequences.

The way in this work can be used to determine the linear complexities and  the trace representations of $r$-ary sequences over $\mathbb{Z}_r$ ($r$ is a prime power) defined by the cyclotomic generator of order $r$  studied in \cite{DH}, the Ding-Helleseth generalized cyclotomic classes of order $r$ modulo $pq$ \cite{CDW}, and the  generalized cyclotomic  classes of order $r$ modulo $p^m$ \cite{DC,YHX2}.

\section*{Acknowledgements}
Parts of this work were written during a very pleasant visit of the
 author to University of Kentucky in Lexington, USA. He wishes to thank Prof. Andrew Klapper for his
 comments and hospitality.

Z.X.C. was partially supported by the National Natural Science
Foundation of China under grant No. 61373140.

\end{document}